\documentclass{article}

\usepackage{amsfonts}
\usepackage{amssymb}
\usepackage[leqno]{amsmath}

\begin{document}

\begin{center}
{\textbf{Vector valued Hardy spaces}} \\\
\end{center}

\begin{center}
Richard D. Carmichael \\
Department of Mathematics \\
Wake Forest University \\
Winston-Salem, N.C. 27109-7388 \\
U.S.A. \\
E-mail: carmicha@wfu.edu \\
\end{center}

\begin{center}
Stevan Pilipovi\'{c} \\
Department of Mathematics and Informatics \\
Faculty of Sciences \\
University of Novi Sad \\
Trg Dositeja Obradovi\'{c}a 4 \\
21000 Novi Sad \\
Serbia \\
E-mail: Stevan.Pilipovic@dmi.uns.ac.rs \\
\end{center}

\begin{center}
Jasson Vindas \\
Department of Mathematics \\
Ghent University \\
Krijgslaan 281 Building S25 \\
B 9000 Gent \\
Belgium \\
E-mail: Jasson.Vindas@ugent.be \\\
\end{center}


{\textit{Acknowledgement}}: The problem considered in this paper was proposed by the first author at the DAAP 2016 (Different Aspects of Analysis and Probability) conference held at the University of Rzesz\'{o}w, Poland, in September 2016. The conference was held to honor the 70th birthday of Professor Dr. Andrzej Kami\'{n}ski. \\\\

{\textbf{Abstract.}} The Hardy space $H^{p}$ of vector valued analytic functions in tube domains in $\mathbb{C}^{n}$ and with values in Banach space are defined. Vector valued analytic functions in tube domains in $\mathbb{C}^{n}$ with values in Hilbert space and which have vector valued tempered distributions as boundary value are proved to be in $H^{p}$ corresponding to Hilbert space if the boundary value is in $L^{p}$ with values in Hilbert space. A Poisson integral representation for such vector valued analytic functions is obtained. \\\\

\noindent {\textbf{1 Introduction}} \\

The representation of tempered distributions as boundary values of analytic functions was first accomplished by H-G Tillmann [22] whose analysis was for functions analytic in half planes and tubes defined by n-rants. Meise ([15], [16]) extended results of this type to vector valued tempered distributions. Carmichael and Walker [2] have obtained a boundary value result for vector valued functions analytic in tubes defined by cones with the boundary value being in the vector valued tempered distribution space. The tempered distributions have also been extended to tempered ultradistributions; see Pilipovi\'{c} [17] and the book [7]. Recently in [10] and [11] a large class of distribution spaces has been introduced and studied whose definition is based on translation-invariant Banach spaces and which generalize the Schwartz $\mathcal{D}'_{L^{p}}$ spaces; a complete theory of boundary value results and analytic representations of these new distributions is obtained. Further, new results associated with the analysis of this paper are contained in [8] and [9]. Boundary value results concerning quasianalytic ultradistributions are obtained in [8], and classical results on boundary values in distribution spaces are important tools in the study of complex Tauberian theorems for Laplace transforms in [9]. The reader should especially note Schwartz ([19], [20]) for a general background of vector valued distributions.

In [18] Raina showed that if the distributional boundary value of an analytic function $f$ in the upper half plane obtained in the dual space of a space of type $\mathcal{S}$ was in fact a $L^{p}$, $1 \leq p \leq \infty$, function then the analytic function $f$ was in the Hardy space $H^{p}$ in the upper half plane. This and associated results found applications in particle physics. In [5] Carmichael and Richters generalized the result of Raina to functions analytic in tube domains and for the tempered and other distribution topologies; associated results and representations were obtained. The main proof of [5] was obtained through representing the assumed analytic function as a Poisson integral in tubes.

As noted above Carmichael and Walker [2, Theorem 8] have obtained a vector valued distributional boundary value result concerning vector valued analytic functions in tubes in $\mathbb{C}^{n}$ which obtain tempered vector valued distributions as boundary value. In this paper we desire to extend the results of [5] to the vector valued case by showing that if the boundary value in [2, Theorem 8] is in fact a $L^{p}$ function with values in a Hilbert space $\mathcal{H}$, then the analytic function, which has values in $\mathcal{H}$, must be in the vector valued Hardy space $H^{p}$ which we define in this paper. \\

\noindent {\textbf{2 Notation and definitions}} \\

Throughout $\mathcal{B}$ will denote a Banach space, $\mathcal{H}$ will denote a Hilbert space, $\mathcal{N}$ will denote the norm of the specified Banach or Hilbert space, and $\Theta$ will denote the zero vector of the specified Banach or Hilbert space. We reference Dunford and Schwartz [12] for integration of vector valued functions and for vector valued analytic functions. For foundational information concerning vector valued distributions we refer to Schwartz ([19], [20]).

The n-dimensional notation used in this paper will be the same as that in [4] and [7]. Basic information concerning cones in $\mathbb{R}^{n}$ can be found in Vladimirov [23] and in the books [6] and [7]. We recall some needed concepts of cones that are important for this paper. $C \subset \mathbb{R}^{n}$ is a cone (with vertex at $\overline{0} = (0,0,...,0)$ in $\mathbb{R}^{n}$) if $y \in C$ implies $\lambda y \in C$ for all positive scalars $\lambda$. The intersection of $C$ with the unit sphere $|y| = 1$ is called the projection of $C$ and is denoted $pr(C)$. A cone $C'$ such that $pr(\overline{C'}) \subset pr(C)$ is a compact subcone of $C$. The function

\begin{equation}
u_{C}(t) = \sup_{y \in pr(C)} (-<t,y>) \, , t \in \mathbb{R}^{n} \, ,
\end{equation}

\noindent is the indicatrix of $C$. The dual cone $C^{*}$ of $C$ is $C^{*} = \{ t \in \mathbb{R}^{n} : \; <t,y> \: \geq \: 0 \; for \; all \; y \in C \} = \{ t \in \mathbb{R}^{n} : u_{C}(t) \leq 0 \}$. A convex cone which does not contain any entire straight line is called a regular cone. Let $v = (v_{1},v_{2},...,v_{n})$ be any of the $2^{n}$ n-tuples whose entries are $0$ or $1$; the $2^{n}$ n-rants $C_{v} = \{ y \in \mathbb{R}^{n}: \, (-1)^{v_{j}}y_{j} > 0, \, j = 1,2,...,n \}$ are examples of regular cones in $\mathbb{R}^{n}$ as are the future and past light cones [23, p. 219]. The n-rants and the light cones are also examples of self dual cones in which the closure of the cone is equal to the dual cone of the cone.

Vector valued functions will be denoted by bold letters; and the space of strongly  measurable  functions  $\bold f$ with the property
$\int_{\mathbb R^n}(\mathcal N(\bold f(t)))^pdt<\infty, p\in[1,\infty),$ is denoted by $L^p(\mathbb R^n,\mathcal B).$ The symbol $|\bold h|_{p}$ will denote the norm of $\textbf{h} \in L^{p}(\mathbb R^n,\mathcal{B})$. In the case $p=\infty$ we make the obvious change of the definition.
Schwartz  spaces $\mathcal S$ and $\mathcal S^{(m)}, m\in\mathbb N,$ over $\mathbb R^n$ are the test spaces; and $\mathcal S'(\mathbb R^n,\mathcal B)$ and
$\mathcal S^{(m)'}(\mathbb R^n,\mathcal B)$ are the spaces of tempered vector  valued distributions with values in $\mathcal B$ with the dual pairing
$$\langle \bold f(t), \phi(t)\rangle\in\mathcal B, \phi\in\mathcal S.
$$

Elements of  $L^p(\mathbb R^n,\mathcal B)$ define the regular elements of $\mathcal S'(\mathbb R^n,\mathcal B)$ by
$$\langle \bold f(t), \phi(t)\rangle= \int_{\mathbb R^n}\bold f(t)\phi(t)dt, \phi\in\mathcal S \; (\mbox{ or }  \phi\in\mathcal S^{(m)}).
$$
We define the Fourier transform in $\mathcal S'(\mathbb R^n,\mathcal B)$ as a transpose mapping
$$\langle \widehat{\bold f}(t), \phi(t)\rangle=\langle \bold f(t), \widehat{\phi}(t)\rangle,\phi\in \mathcal S,
$$

\noindent where the Fourier transform of $\phi \in \mathcal{S}$ (or for any $L^{1}(\mathbb{R}^{n})$ function) is

\begin{center}
$\widehat{\phi}(x) = \mathcal{F}({\phi})(x) = \mathcal{F}[\phi (t);x] = \int_{\mathbb{R}^{n}} e^{2 \pi i<x,t>} \phi (t)dt, \; x \in \mathbb{R}^{n}.$
\end{center}

Appropriate references for vector valued weighted functions in this paper are [1] and [14]. If $\bold f\in L^1(\mathbb R^n,\mathcal B)$, then its Fourier transform is defined by  the Bochner integral
$$\mathcal F(\bold f)(x)=\mathcal{F}[{\bold f}(t);x] =\int_{\mathbb R^n}e^{2\pi i <x,t>}\bold f(t)dt, x\in\mathbb R^n ;$$
\noindent while $\mathcal F^{-1}(\bold f)(x)=\mathcal F(\bold f)(-x)$.
Then ([1], [14]) $\mathcal F(\bold f)\in L^\infty(\mathbb R^n,\mathcal B)$ and the Hausdorff-Young inequality holds. Moreover, $\widehat{\bold f}\in C_0(\mathbb R^n,\mathcal B)$
which means that $\widehat{\bold f}$ is continuous and $\mathcal N(\widehat{\bold f}(x))\rightarrow 0$ as $|x|\rightarrow \infty.$

For the case $p=2$ in order to have an isomorphism of $\mathcal F $ of $ L^2(\mathbb R^n,\mathcal B) $ onto itself
(with the  Parseval identity  $||\mathcal N(\widehat{\bold f})||_{L^2(\mathbb R^n)}=|| \mathcal N(\bold f)||_{L^2(\mathbb R^n)}$)
 it is necessary and sufficient that $\mathcal B=\mathcal H$ is a Hilbert space. Recall that in this case $\mathcal F$ is defined first on a dense set  and then extended on $ L^2(\mathbb R^n,\mathcal H) $ [1].

By the Fubini theorem ([1]),  the Fourier transform of $\bold f$ in the cases  of $L^1(\mathbb R^n,\mathcal B)$ and $L^2(\mathbb R^n, \mathcal H)$  defined in the sense of Bochner integral imbedded into the
$\mathcal S'(\mathbb R^n,\mathcal B)$ coincide with the $\widehat{\bold f}$ defined as the distributional Fourier transform.
We will use this fact in the sequel.

The Hardy space $H^{p}(T^{C},\mathcal{B})$, $1 \leq p < \infty$, consists of those analytic functions $\bold f(z)$ in the tube $T^{C} = \mathbb{R}^{n} + iC$ with values in $\mathcal{B}$  such that

\begin{center}
$\int_{\mathbb{R}^{n}} (\mathcal{N}(\bold f(x+iy)))^{p}dx \; \leq \; A, \, y \in C,$
\end{center}

\noindent where the constant $A$ is independent of $y \in C$. The Hardy space $H^{\infty}(T^{C},\mathcal{B})$  is defined with the usual changes.

Let $C$ be an open convex cone in $\mathbb{R}^{n}$. The following function $d_{y}(t), t \in \mathbb{R}^{n}, y \in C$, will be used. Let $s(u) \in \mathcal{E}(\mathbb{R}), u \in \mathbb{R}$, such that $s(u) = 1, u \geq 0, s(u) = 0, u \leq -{\epsilon}, {\epsilon} > 0$ and fixed, and $0 \leq s(u) \leq 1$. Put

\begin{equation}
d_{y}(t) = s(<t,y>), \; t \in \mathbb{R}^{n}, \; y \in C.
\end{equation}

\noindent We have $d_{y}(t) \in \mathcal{E}(\mathbb{R}^{n})$, for any $y \in C$. \\

\noindent {\textbf{3 Preliminary results}}  \\

We indicate results in this section which we need to prove the main results of this paper contained in section 4. Recall that $C \subset \mathbb{R}^{n}$ is a regular cone if it is an open convex cone which does not contain any entire straight line.

For $C$ being a regular cone, the Cauchy kernel corresponding to the tube $T^{C} = \mathbb{R}^{n} + iC$ is

\begin{center}
$K(z-t) = \int_{C^{*}} e^{2 \pi i<z-t, \eta>} d \eta, \, t \in \mathbb{R}^{n}, \, z \in T^{C},$
\end{center}

\noindent where $C^{*}$ is the dual cone of $C$. The Poisson kernel corresponding to $T^{C}$ is

\begin{center}
$Q(z;t) = \frac{K(z-t) \overline{K(z-t)}}{K(2iy)} = \frac{|K(z-t)|^{2}}{K(2iy)}, \, t \in \mathbb{R}^{n}, \, z \in T^{C}.$
\end{center}

Referring to [7] for details we know that $K(z-\cdot) \in \mathcal{D}(*,L^{p})\subset\mathcal D_{L^p}, 1 < p \leq \infty$, and $Q(z;\cdot) \in \mathcal{D}(*,L^{p})\subset\mathcal D_{L^p}, (z \in T^{C}) 1 \leq p \leq \infty$, where $*$ is Beurling $(M_{p})$ or Roumieu $\{ M_{p} \}$. These ultradifferentiable functions are contained in the Schwartz space $\mathcal{D}(L^{p},\mathbb{R}^{n})$. We will use the results of [5, Lemmas 3.1 and 3.2] here as well as the calculation [5, (4.8)] and calculations in [5, section 2], and they will not be restated here.

Because of these properties of the Cauchy and Poisson kernels, we know that the Cauchy and Poisson integrals

\begin{center}
$\int_{\mathbb{R}^{n}} \bold h(t)K(z-t)dt \; and \; \int_{\mathbb{R}^{n}} \bold h(t)Q(z;t)dt, \, z \in T^{C},$
\end{center}

\noindent are well defined for $\bold h \in L^{p}(\mathbb R^n,\mathcal{B}), 1 \leq p < \infty$, and $\bold h \in L^{p}(\mathbb R^n,\mathcal{B}),
1 \leq p \leq \infty$, respectively.

Several lemmas are proved now which are needed for section 4. \\

LEMMA 3.1. {\textit{Let $C$ be an open connected cone and $C'$ be an arbitrary compact subcone of $C$. Let $r > 0$ be arbitrary.
Let $\bold g(t), t \in \mathbb{R}^{n}$, be a continuous function with support in $C^{*}$ and with values in a Banach space $\mathcal{B}$ which satisfies}}

\begin{equation}
\mathcal{N}(\bold g(t)) \leq M(C',r) exp(2 \pi (<y',t>+ \sigma |y'|)), \, t \in \mathbb{R}^{n},
\end{equation}

\noindent {\textit{for all $\sigma > 0$, where $M(C',r)$ is a constant which depends on $C'$ and on $r > 0$ and (3) is independent of $y' \in (C' \setminus (C' \cap N(\overline{0},r)))$ (that is, (3) holds for all $y' \in (C' \setminus (C' \cap N(\overline{0},r))))$. Let $y$ be any point of $C$. We have}}
$$(e^{-2 \pi <y,\cdot>}\bold g) \in L^{p}(\mathbb R^n,\mathcal{B})$$ {\textit{for all
$p, 1 \leq p < \infty$.}} \\

{\textit{Proof}}. Let $y \in C$ be arbitrary but fixed. There exists $C' \subset C$ and $r > 0$ such that $y \in (C' \setminus (C' \cap N(\overline{0},r)))$. Choose $\mu$ such that $1 > \mu >(r/|y|) > 0$ and put $y' = \mu y$. We have $y' \in C'$ since $C'$ is a cone and $|\mu y| = \mu |y| > r > 0$; thus $ y' \in C' \setminus (C' \cap N(\overline{0},r))$. By [6, Lemma 4.3.2, p. 155] there is a $\delta > 0$ such that

\begin{equation}
<y,t> \, \geq \, \delta |y||t|, \; t \in C^{*}, \; y \in C',
\end{equation}

\noindent and $\delta$ depends only on $C'$ and not on $y \in C'$. Now taking $y' = \mu y$ in (3) we have for $t \in C^{*}$ and $y \in C'$

\begin{eqnarray}
 &~&\mathcal{N}(e^{-2 \pi <y,t>}\bold g(t))  = e^{-2 \pi <y,t>} \mathcal{N}(\bold g(t)) \nonumber \\
 &~&  \leq M(C',r) exp(2 \pi (\mu <y,t>+ \sigma \mu |y|-<y,t>)) \nonumber \\
 &~&  = M(C',r) e^{2 \pi \sigma \mu |y|}exp(2 \pi (1- \mu )(-<y,t>)) \nonumber \\
 &~& \leq M(C',r) e^{2 \pi \sigma \mu |y|} exp(-2 \pi \delta (1- \mu)|y||t|) \nonumber
\end{eqnarray}

\noindent with $\delta > 0$ and $(1- \mu) > 0$. Since $supp(\bold g) \subseteq C^{*}$, we have for $1 \leq p < \infty$ and the arbitrary but fixed $y \in C$ at the beginning of the proof

\begin{eqnarray}
&~& \int_{\mathbb{R}^{n}} (\mathcal{N}(e^{-2 \pi <y,t>}\bold g(t)))^{p}dt = \int_{C^{*}} (\mathcal{N}(e^{-2 \pi <y,t>}\bold g(t)))^{p}dt \\ &~& \leq (M(C',r))^{p} e^{2 \pi \sigma \mu p|y|} \int_{C^{*}} exp(-2 \pi \delta (1- \mu)p|y||t|) dt \nonumber \\ &~& \leq Z_{n} (M(C',r))^{p} e^{2 \pi \sigma \mu p|y|} \int_{0}^{\infty} u^{n-1} exp(-2 \pi \delta (1- \mu)p|y|u) du \nonumber \\ &~& = Z_{n} (M(C',r))^{p} e^{2 \pi \sigma \mu p|y|} (n-1)! (2 \pi \delta (1- \mu)p|y|)^{-n} \nonumber
\end{eqnarray}

\noindent where $Z_{n}$ is the surface area of the unit sphere in $\mathbb{R}^{n}$, which proves $e^{-2 \pi <y,t>}\bold g(t) \in L^{p}(\mathbb R^n,\mathcal{B}), 1 \leq p < \infty,$ for any $y \in C$. \\

In the following Lemmas 3.2, 3.3, and 3.4 $C$ will be a regular cone. The vector valued functions in Lemma 3.2 have values in a Hilbert space as opposed to a general Banach space because we use the Fourier transform in the proof; recall that if $\mathcal B = \mathcal H$, a Hilbert space, the Fourier transform is a bijection on $L^{2}(\mathbb{R}^{n},\mathcal H)$. For Lemma 3.4 $\mathcal B$ is an arbitrary Banach space. \\

LEMMA 3.2. {\textit{Let $\bold g \in L^{2}(\mathbb R^n,\mathcal{H})$ and $\bold G(\eta) = \mathcal{F}^{-1}[\bold g(t); \eta], \eta\in\mathbb R^n,$ in $L^{2}(\mathbb R^n,\mathcal{H})$. Assume $\bold Gexp(2 \pi i<z,\cdot>) \in L^{1}(\mathbb R^n,\mathcal{H})$ for $z \in T^{C}$ and that $supp(\bold G) \subseteq C^{*}$ almost everywhere. We have}}

\begin{center}
$\int_{\mathbb{R}^{n}} \bold G(\eta) e^{2 \pi i<z,\eta>}d \eta = \int_{\mathbb{R}^{n}} \bold g(t) K(z-t) dt , \; z \in T^{C}.$ \\
\end{center}

{\textit{Proof.}} Let $I_{C^{*}}(\eta)$ denote the characteristic function of $C^{*}$. By [5, Lemma 2.1] $I_{C^{*}}(\eta)e^{2 \pi i<z,\eta>} \in L^{p}, 1 \leq p \leq \infty$, as a function of $\eta \in \mathbb{R}^{n}$ for $z \in T^{C}$. Recalling the Cauchy kernel $K(z-t), t \in \mathbb{R}^{n}, z \in T^{C}$, we have

\begin{eqnarray}
&~& \int_{\mathbb{R}^{n}} \bold g(t) K(z-t) dt = \int_{\mathbb{R}^{n}} \bold g(t) \int_{C^{*}} e^{2 \pi i<z-t,\eta>} d \eta dt \nonumber \\ &~& = \int_{\mathbb{R}^{n}} \bold g(t) \mathcal{F}^{-1}[I_{C^{*}}(\eta) e^{2 \pi i<z,\eta>};t] dt = \int_{C^{*}} \mathcal{F}^{-1}[\bold g(t);\eta] e^{2 \pi i<z,\eta>} d \eta \nonumber \\ &~& = \int_{C^{*}} \bold G(\eta) e^{2 \pi i<z,\eta>} d \eta = \int_{\mathbb{R}^{n}} \bold G(\eta) e^{2 \pi i<z, \eta>} d \eta . \nonumber
\end{eqnarray}

LEMMA 3.3. {\textit{Let $z_{o}$ be an arbitrary but fixed point in $T^{C}$. Let $1 \leq p \leq \infty$. There exists a closed neighborhood $N(z_{o}, \rho) = \{ z: |z-z_{o}| \leq \rho, \; \rho > 0 \}$ of $z_{o}$ which is contained in $T^{C}$ and a constant $B(z_{o})$ depending only on $z_{o}$ such that}}

\begin{center}
$||Q(z;t)||_{L^{p}} \; \leq \; B(z_{o}) \; < \; \infty, \; z \in N(z_{o},\rho),$
\end{center}

\noindent {\textit{where the $L^{p}$ norm is with respect to $t \in \mathbb{R}^{n}$ and $Q(z;t), t \in \mathbb{R}^{n}, \, z \in T^{C},$ is the Poisson kernel.}} \\

{\textit{Proof.}} See [5, Lemma 3.4]. \\

LEMMA 3.4. {\textit{Let $\bold f$ be analytic in $T^{C}$ with values in a Banach space $\mathcal{B}$ ($\bold f \in \mathcal{A}(T^{C},\mathcal{B}))$ and have the Poisson integral representation}}

\begin{center}
$\bold f(z) = \int_{\mathbb{R}^{n}} \bold h(t) Q(z;t) dt, \; z \in T^{C},$
\end{center}

\noindent {\textit{for ${\textbf{h}}  \in L^{p}(\mathbb R^n,\mathcal{B}), 1 \leq p \leq \infty$. We have $\bold f \in H^{p}(T^{C},\mathcal{B}), 1 \leq p \leq \infty$. For $p = \infty, \; \bold f(x+iy) \rightarrow \bold h(x)$ in the weak-star topology of $L^{\infty}(\mathbb{R}^{n},\mathcal{B})$ as $y \rightarrow \overline{0}, y \in C; \; for \; 1 \leq p < \infty, \; \bold f(x+iy) \rightarrow \bold h(x), x \in \mathbb{R}^n$, as $y \rightarrow \overline{0}, y \in C,$ in $L^{p}(\mathbb{R}^n,\mathcal{B})$; \; for $1 < p \leq 2$}}

\begin{equation}
\mathcal{N}(\bold f(x+iy)) \leq M(C') |\bold h|_{p}|y|^{-n/p}, \; z = x+iy \in T^{C'},
\end{equation}

\noindent {\textit{for all compact subcones $C' \subset C, \; M(C')$ being a constant depending on $C' \subset C$ and not on $y \in C'$, while}}

\begin{center}
$\mathcal{N}(\bold f(x+iy)) \leq M_{y} |\bold h|_{p}|y|^{-n/p}, \; z = x+iy \in T^{C},$
\end{center}

\noindent{\textit{ where $M_{y}$ is a constant depending on $y \in C$; and for $2 < p < \infty$}}

\begin{eqnarray}
&~&\mathcal{N}(\bold f(x+iy)) \leq M(C',r) |\bold h|_{p}, \nonumber \\ &~& z = x+iy \in T(C',r) \; = \; \{ z = x+iy: x \in \mathbb{R}^{n}, \; y \in (C' \setminus (C' \cap N(\overline{0},r))) \}, \nonumber
\end{eqnarray}

\noindent {\textit{for all compact subcones $C' \subset C$ and all $r > 0, M(C',r)$ being a constant depending on $C' \subset C$ and on $r > 0$ but not on $y \in (C' \setminus (C' \cap N(\overline{0},r)))$ while}}

\begin{center}
$ \mathcal{N}(\bold f(x+iy)) \leq M_{y} |\bold h|_{p}, \; z = x+iy \in T^{C},$
\end{center}

\noindent {\textit{where $M_{y}$ is a constant depending on $y \in C$.}} \\

{\textit{Proof.}} For $p = \infty$ and $z \in T^{C}$ we use [5, (3.3) and (3.4)] to obtain

\begin{center}
$\mathcal{N}(\bold f(z)) \leq \int_{\mathbb{R}^{n}} \mathcal{N}(\bold h(t)) Q(z;t)dt \; \leq \; A$
\end{center}

\noindent where $A$ is the bound on $\bold h$; and $\bold f \in H^{\infty}(T^{C},\mathcal{B}).$ Also for $p = \infty$ the weak-star convergence of $\bold f(x+iy)$ to $\bold h(x)$ as $y \rightarrow \overline{0}, \; y \in C,$ is proved as in the scalar valued case using the approximate identity properties of the Poisson kernel; see [5, Lemma 3.5]. For $1 \leq p < \infty$ we use Jensen's inequality [13, 2.4.19, p. 91], which holds for Banach spaces $\mathcal{B}$, and [5, Lemma 3.1] to obtain for $y \in C$

\begin{eqnarray}
&~& \int_{\mathbb{R}^{n}} (\mathcal{N}(\bold f(x+iy)))^{p}dx \leq \int_{\mathcal{R}^{n}} \int_{\mathbb{R}^{n}} (\mathcal{N}(\bold h(t)))^{p} Q(z;t) dt dx \nonumber \\ &~& = \int_{\mathbb{R}^{n}} (\mathcal{N}(\bold h(t)))^{p} \int_{\mathbb{R}^{n}} Q(z;t) dx dt = \int_{\mathbb{R}^{n}} (\mathcal{N}(\bold h(t)))^{p} dt \; < \, \infty, \nonumber
\end{eqnarray}

\noindent and $\bold f \in H^{p}(T^{C},\mathcal{B})$.

We use the properties of $Q(z;t)$ in [5, Lemma 3.1] and the associated properties of $Q(u;y), u \in \mathbb{R}^{n}, y \in C$, stated in [21, p. 105]. For $z = x+iy \in T^{C}$ and $\rho > 0$

\begin{eqnarray}
&~& |\bold f(x+iy) - \bold h(x)|_{p} = \left( \int_{\mathbb{R}^{n}} (\mathcal{N} ( \int_{\mathbb{R}^{n}} \bold h(t)Q(z;t)dt - \int_{\mathbb{R}^{n}} \bold h(x)Q(z;t)dt))^{p}dx \right)^{1/p} \nonumber \\  &~& \leq ( \int_{\mathbb{R}^{n}} 2^{p} ((\mathcal{N} ( \int_{|u| \leq \rho} (\bold h(x-u)-\bold h(x))Q(u;y)du))^{p} \nonumber \\ &~& + ( \mathcal{N} ( \int_{|u| > \rho} (\bold h(x-u)-\bold h(x))Q(u;y)du))^{p})dx )^{1/p} \nonumber \\  &~&  \leq 2^{(p+1)/p} ( \int_{|u| \leq \rho} \left(\int_{\mathbb{R}^{n}} (\mathcal{N} (\bold h(x-u)-\bold h(x)))^{p} dx \right)^{1/p} Q(u;y)du \nonumber \\ &~& + \int_{|u| > \rho} \left( \int_{\mathbb{R}^{n}} ( \mathcal{N} (\bold h(x-u)-\bold h(x)))^{p}dx \right)^{1/p} Q(u;y)du) \nonumber \\ &~& \leq 2^{(p+1)/p} \left( \sup_{|u| \leq \rho} \left( \int_{\mathbb{R}^{n}} ( \mathcal{N}(\bold h(x-u)-\bold h(x)))^{p} dx \right)^{1/p} \right) \: \int_{\mathbb{R}^{n}} Q(u;y)du \nonumber \\ &~& + 2^{(3p+2)/p} |\bold h|_{p} \int_{|u| > \rho} Q(u;y)du \; . \nonumber
\end{eqnarray}

\noindent The sup term $\rightarrow 0$ as $|u| \rightarrow 0$; thus we can choose $\rho > 0$ such that

\begin{center}
$2^{(p+1)/p} ( \sup_{|u| \leq \rho} ( \int_{\mathbb{R}^{n}} ( \mathcal{N} (\bold h(x-u)-\bold h(x)))^{p} dx)^{1/p}) \; < \; \epsilon \; $
\end{center}

\noindent for $\epsilon > 0$. For this chosen $\rho$ we have

\begin{center}
$ \int_{|u| > \rho} Q(u;y)du \; \rightarrow \; 0$
\end{center}

\noindent as $y \rightarrow \overline{0}, y \in C,$ from [21, p. 105]. Combining we obtain $\bold f(x+iy) \rightarrow \bold h(x)$ in $L^{p}(\mathbb{R}^{n},\mathcal{B})$ as $y \rightarrow \overline{0}, y \in C$.

For the remainder of the conclusions recall from section $3$ that for $z \in T^{C}, K(z-t) \in \mathcal{D}(*,L^{p}) \subset \mathcal{D}(L^{p},\mathbb{R}^{n}), 1 < p \leq \infty$, and $Q(z;t) \in \mathcal{D}(*,L^{p}) \subset \mathcal{D}(L^{p},\mathbb{R}^{n}), 1 \leq p \leq \infty$, where $*$ is either Beurling or Roumieu. Now let $1 < p \leq 2$. By H$\ddot{o}$lder's inequality

\begin{equation}
\mathcal{N} (\bold f(x+iy)) \; \leq \; |\bold h|_{p} ||Q(x+iy;t)||_{L^{q}}, \; z = x+iy \in T^{C},
\end{equation}

\noindent $1/p + 1/q = 1$. From the definition of the Poisson kernel for $z \in T^{C}$

\begin{equation}
||Q(x+iy;t)||_{L^{q}} = (K(2iy))^{-1} ( \int_{\mathbb{R}^{n}} |K(x+iy-t)|^{2q} dt)^{1/q}, \; z = x+iy \; \in \; T^{C}.
\end{equation}

\noindent By [5, Lemma 2.1] $I_{C^{*}}(\eta) e^{2 \pi i<z,\eta>} \; \in \; L^{1} \cap L^{p}, \; 1 < p \: \leq 2$, for $z \in T^{C}$. Thus

\begin{center}
$K(z-t) \; = \; \mathcal{F}^{-1}[I_{C^{*}}(\eta)e^{2 \pi i<z,\eta>};t], \; z \; \in \; T^{C},$
\end{center}

\noindent and by the Parseval inequality

\begin{center}
$||K(z-t)||_{L^{q}} \; \leq \; ||I_{C^{*}}(\eta) e^{2 \pi i<z,\eta>}||_{L^{p}} \, , \; z \; \in \; T^{C}.$
\end{center}

\noindent By this Parseval inequality and analysis as in (5) we have

\begin{equation}
||K(x+iy-t||_{L^{q}} \; \leq \; \left( \frac{Z_{n}(n-1)!}{(2 \pi p \delta)^{n}} \right)^{1/p} \; |y|^{-n/p}
\end{equation}

\noindent with $\delta$ depending on $y \in C$ if (9) is taken to hold for all $y \in C$ while $\delta$ depends on $C' \subset C$, and not on $y \in C' \subset C$, if (9) is taken to hold for $y \in C' \subset C$ for $C'$ being any compact subcone of $C$. From [4, Lemma 2],

\begin{equation}
K(2iy) \; \geq \; B(C) |y|^{-n}, \; y \; \in \; C,
\end{equation}

\noindent where the constant $B(C)$ depends only on $C$ and not on $y \in C$. Additionally from [3, Lemma 3],

\begin{equation}
|K(z-t)| \; \leq \; Z_{n} (n-l)! \delta^{-n} |y|^{-n} \, , \; t \; \in \; \mathbb{R}^{n},
\end{equation}

\noindent with $\delta > 0$ depending on $y \in C$ and (11) holding for $z = x+iy \in T^{C}$ while (11) holds for all $z = x+iy \in T^{C'}, C'$ being an arbitrary compact subcone of $C$, with $\delta$ depending only on $C' \subset C$ and not on $y \in C' \subset C$. Now the desired norm growth on $\mathcal{N}(\bold f(\cdot+iy))$ for $1 < p \leq 2$ follows by combining (7), (8), (9), (10), and (11).

Now let $2 < p < \infty$. Again by H$\ddot{o}$lder's inequality as in (7), we combine (8), (10), and (11) to obtain

\begin{center}
$\mathcal{N}(\bold f(\cdot+iy)) \leq |\bold h|_{p} (|y|^{n}/B(C))(Z_{n}(n-1)!/\delta^{n}|y|^{n})||K(\cdot+iy-t)||_{L^{q}}$
\end{center}

\noindent with the $\delta > 0$ depending on $C' \subset C$ for $y \in C' \subset C$ and not on $y$ while $\delta > 0$ depends on $y$ if $y \in C$. From the proof of [7, Theorem 4.1.1] or by using (9), for $z = x+iy \in T(C',r) = \{ z = x+iy: x \in \mathbb{R}^{n}, \; y \in (C' \setminus (C' \cap N(\overline{0},r))) \}$ where $C'$ is any compact subcone of $C$ and $r > 0$ is arbitrary we have

\begin{center}
$||K(x+iy-t)||_{L^{q}} \leq P(C',r), \; z = x+iy \in T(C',r),$
\end{center}

\noindent while

\begin{center}
$||K(x+iy-t)||_{L^{q}} \leq P_{y}, \; z = x+iy \in T^{C},$
\end{center}

\noindent with the constant $P_{y}$ depending on $y \in C$. Combining the inequalities we have the desired estimates for $\mathcal{N}(\bold f(x+iy))$ when $2 < p < \infty$. \\

The following lemma is proved by the same proof as [5, Lemma 3.6]. \\

LEMMA 3.5. {\textit{Let $\bold h \in L^{\infty}(\mathbb{R}^{n}, \mathcal{B})$. Let $C$ be a regular cone. Put}}

\begin{center}
$X_{\epsilon}(t) = \prod_{j=1}^{n} (1-i \epsilon (-1)^{v_{j}}t_{j})^{R+n+2}, \; \epsilon > 0, \; t \in \mathbb{R}^{n},$
\end{center}

\noindent {\textit{where $R \geq 0$ is a fixed real number, $n$ is the dimension, and $v = (v_{1},v_{2},...,v_{n})$ is any of the $2^{n}$ n-tuples whose entries are $0$ or $1$ that defines the quadrant $C_{v}$. We have}}

\begin{center}
$\lim_{\epsilon \rightarrow 0+} \mathcal{N}(\int_{\mathbb{R}^{n}} (\bold h(t) - \frac{\bold h(t)}{X_{\epsilon}(t)})Q(z;t)dt) \; = \; 0$
\end{center}

\noindent {\textit{uniformly in $z$ on compact subsets of $T^{C}$.}} \\

\noindent {\textbf{4 $H^{p}(T^{C},\mathcal{H})$ functions, $2 \leq p \leq \infty$}} \\

We begin by restating [2, Theorem 8, p. 327] in a more general form which is proved by the same analysis of [2, Theorem 8, p. 327]. In Theorem 4.1 $\mathcal{B}$ can be any space assumed in [2, Theorem 8]. However $\mathcal{B}$ being a Banach space is of primary interest in this paper. \\

THEOREM 4.1. {\textit{Let $C$ be an open convex cone. Let $\bold f \in \mathcal A(T^{C}, \mathcal B)$. For every compact subcone $C' \subset C$ and every $r > 0$ let}}

\begin{eqnarray}
&~& \mathcal{N}(\bold f(x+iy)) \leq M(C',r) (1+|x|)^{R}|y|^{-k}, \\ &~& z = x+iy \; \in \; T(C',r) \; = \; \mathbb{R}^{n} + i(C' \setminus (C' \cap N(\overline{0},r))), \nonumber
\end{eqnarray}

\noindent {\textit{where $M(C',r)$ is a constant depending on $C' \subset C$ and on $r$, $R$ is a nonnegative integer, $k$ is an integer greater than $1$, and neither $R$ nor $k$ depend on $C'$ or $r$. There exists a positive integer $m$ and a unique element $\bold U \in \mathcal{S}^{(m)'}(\mathbb R^n,\mathcal{B}) \subset \mathcal{S}'(\mathbb R^n,\mathcal{B})$ such that}}

\begin{equation}
\lim_{y \rightarrow \overline{0}, y \in C} \mathcal{N}(\langle \bold f(x+iy), \phi (x) \rangle \; - \; \langle \bold U, \phi \rangle ) \; = \; 0, \; \phi \in \mathcal{S}^{(m)}. \\
\end{equation}

$\bold U$ here will be called the $\mathcal{S}^{(m)'}(\mathbb R^n,\mathcal{B}) \subset \mathcal{S}'(\mathbb R^n,\mathcal{B})$ boundary value of $\bold f(\cdot+iy)$. In Theorem 4.1 and in the other theorems in this section, by $y \rightarrow \overline{0}, \; y \in C$, we mean that $y \rightarrow \overline{0}, \; y \in C' \subset C$ for every compact subcone $C'$ of $C$.

The element $\bold U$ in the conclusion of Theorem 4.1 could be a function $\bold h \in L^{p}(\mathbb R^n,\mathcal{B})$.
We now prove that if this is the case for $\mathcal{B}$ being a Hilbert space $\mathcal{H}$ and $2 \, \leq \, p \, \leq \, \infty$, the analytic function $\bold f$ in Theorem 4.1 is in fact an element of $H^{p}(T^{C},\mathcal{H})$. We prove this in two steps. First we consider the cone $C$ to be contained in or be any of the $2^{n}$ n-rants $C_{v}$; we then use this case to prove our result for $C$ being any regular cone.

Before proceeding to the statement and proof of our desired result for the case that the cone $C$ satisfies $C \subseteq C_{v}$, we first give an outline of the proof. First note that the functions and distributions are assumed to have values in a Hilbert space $\mathcal{H}$ now as opposed to a general Banach space; the reason for this is that we use the function Fourier transform in the proof. We first consider the case $p = 2$. Starting with the assumed analytic function $\bold f(z), \; z \in T^{C}$, which has boundary value $\bold h \in L^{2}(\mathbb R^n,\mathcal{H})$ in $\mathcal{S}'(\mathbb R^n,\mathcal{H})$, we multiply $\bold f(z)$ by an analytic function depending on $\epsilon > 0$ such that the product $\bold g_{\epsilon}(z)$ has a stronger growth than $\bold f(z)$ and has a boundary value $\bold U_{\epsilon} \; \in \; \mathcal{S}'(\mathbb R^n,\mathcal{H})$ as $y \; \rightarrow \; \overline{0}, \; y \in C, \; \epsilon > 0$. A function of $t \in \mathbb{R}^{n}$ depending on $\epsilon > 0, \; \bold G_{\epsilon}(t)$, is constructed as the Fourier-Laplace transform of $\bold g_{\epsilon}(z)$ and needed properties of $\bold G_{\epsilon}(t)$ are obtained. We show $\bold U_{\epsilon} \; = \; \mathcal{F}[\bold G_{\epsilon}]$ in $\mathcal{S}'(\mathbb R^n,\mathcal{H})$ and proceed to construct a function $\bold H_{\epsilon} \; \in \; L^{2}(\mathbb R^n,\mathcal{H})$, which equals $\bold G_{\epsilon}$ in $\mathcal{S}'(\mathbb R^n,\mathcal{H})$, from which we show that the Poisson integral of $\mathcal{F}[\bold H_{\epsilon}(x);t]$ in $L^{2}$ equals $\bold g_{\epsilon}(z), \; z \; \in \; T^{C}$. From this Poisson integral representation of $\bold g_{\epsilon}(z)$ and its analyticity in $T^{C}$ we prove that the Poisson integral of the assumed boundary value $\bold h \; \in \; L^{2}(\mathbb R^n,\mathcal{H})$ in the theorem is analytic in $T^{C}$ using a limit argument and that this limit Poisson integral is in $H^{2}(T^{C}, \mathcal{H})$. To conclude the proof we show that the original assumed $\bold f(z)$ equals this limit Poisson integral of $\bold h \; \in \; L^{2}(\mathbb R^n,\mathcal{H})$ for $z \; \in \; T^{C}$ and is thus in $H^{2}(T^{C},\mathcal{H})$. The proof of our result for $2 < p \leq \infty$ follows by the same procedure as for the case $p = 2$. \\

THEOREM 4.2. {\textit{Let $C$ be an open convex cone which is contained in or is any of the $2^{n}$ n-rants $C_{v}$ in $\mathbb{R}^{n}$. Let $f \in \mathcal A(T^{C},\mathcal H)$  which satisfies (12). Let the unique boundary value $\bold U$ of Theorem 4.1 be $\bold h \; \in \; L^{2}(\mathbb R^n,\mathcal{H})$. We have $\bold f \; \in \; H^{2}(T^{C},\mathcal{H})$ and}}

\begin{equation}
\bold f(z) \; = \; \int_{\mathbb{R}^{n}}\bold h(t) Q(z;t) dt, \; z \; \in \; T^{C}. \\
\end{equation}

{\textit{Proof.}} Put $\bold g_{\epsilon}(z) \, = \, \bold f(z)/X_{\epsilon}(z), \, z \in T^{C}, \, \epsilon > 0$, where

\begin{center}
$X_{\epsilon}(z) \; = \; \prod_{j=1}^{n} (1 - i \epsilon (-1)^{v_{j}} z_{j})^{R+n+2}, \; \epsilon > 0.$
\end{center}

\noindent  Clearly, $\bold g_{\epsilon}$ satisfies (12). By Theorem 4.1 there is a unique $\bold U_{\epsilon} \in \mathcal{S}'(\mathbb R^n,\mathcal{H})$ such that $\bold g_{\epsilon}(x+iy) \; \rightarrow \; \bold U_{\epsilon}(x), x\in\mathbb R^n,$ in $\mathcal{S}'(\mathbb R^n,\mathcal{H})$ as $y \, \rightarrow \, \overline{0}, \,y \, \in \, C$; that is, (13) holds for $\bold g_{\epsilon}$ and $\bold U_{\epsilon}$. From (12) and the calculations of [5, (4.8)] there is a constant $M'(C',r,\epsilon)$ such that

\begin{eqnarray}
&~& \mathcal{N}(\bold g_{\epsilon}(z)) \leq M'(C',r,\epsilon) (1+|z|)^{-n-2}, \\ &~& z \; \in \; T(C',r) = \mathbb{R}^{n}+i(C' \setminus (C' \cap N(\overline{0},r))), \nonumber
\end{eqnarray}

\noindent for all compact subcones $C' \subset C$ and all $r > 0$. Put

\begin{equation}
\bold G_{\epsilon}(t) = \int_{\mathbb{R}^{n}} \bold g_{\epsilon}(x+iy) e^{-2 \pi i<x+iy,t>}dx, \; y \in C, \; t \in \mathbb{R}^{n}.
\end{equation}

\noindent For any $y \, \in \, C, \; y \, \in \, C' \subset C$ and $|y| > r$ for some compact subcone $C' \subset C$ and some $r > 0$; thus $\bold G_{\epsilon}(t)$ is a well defined function of $t \in \mathbb{R}^{n}$ for any $y \, \in \, C$ and any $\epsilon > 0$ and is a continuous function of $t \, \in \, \mathbb{R}^{n}$ for $y \, \in \, C$ and $\epsilon > 0$. Let $C''$ be an arbitrary compact subdomain of $C$. From (15),

\begin{center}
$\int_{C''} \mathcal{N}(\bold g_{\epsilon}(x+iy) e^{-2 \pi i<x+iy,t>})dy \; \rightarrow \; 0$
\end{center}

\noindent as $|x| \, \rightarrow \, \infty$; hence an application of the Cauchy-Poincare theorem yields that $\bold G_{\epsilon}$ is independent of $y \in C''$ and is thus independent of $y \in C$ since $C''$ is an arbitrary compact subdomain of $C$.

We now show that $supp(\bold G_{\epsilon}) \, \subseteq \, C^{*} \, = \, \{ t \in \mathbb{R}^{n}: <t,y> \; \geq \; 0 \; for \; all \; y \in C \} \; = \; \{ t \in \mathbb{R}^{n}: u_{C}(t) \; \leq \; 0 \}$. Let $t_{o} \; \in \; \mathbb{R}^{n} \setminus C^{*} \; = \; C_{*}$; thus $u_{C}(t_{o}) \; > \; 0$. By the proof of [23, Lemma, p. 241] there is a point $y' \; \in \; pr(C)$ and a number $\rho \; = \; \rho(t_{o}) \; > \; 0$ which can be chosen small enough in order that

\begin{center}
$-<t_{o},y'> \; \geq \; u_{C}(t_{o}) - \rho \; > \; 0$.
\end{center}

\noindent Letting $\lambda \, > \, 0$ be arbitrary, set $C' \, = \, \{ y: y = \lambda y' \} \, \subset \, C$ and recall that $\bold G_{\epsilon}$ is independent of $y \, \in \, C$. Using (15) we have for $\lambda \, > \, r \, > \, 0$

\begin{eqnarray}
&~& \mathcal{N}(\bold G_{\epsilon}(t_{o})) \leq \int_{\mathbb{R}^{n}} \mathcal{N}(\bold g_{\epsilon}(x+i \lambda y')e^{-2 \pi i<x+i \lambda y',t_{o}>})dx \\ &~& \leq M'(C',r,\epsilon) e^{2 \pi \lambda <y',t_{o}>} \int_{\mathbb{R}^{n}} (1+|x|)^{-n-2}dx \nonumber \\ &~& \leq M''(C',r,\epsilon) e^{2 \pi \lambda (\rho - u_{C}(t_{o}))} \nonumber
\end{eqnarray}

\noindent with $(\rho - u_{C}(t_{o})) \, < \, 0$. Letting $\lambda \; \rightarrow \; \infty$ (17) yields $\mathcal{N}(\bold G_{\epsilon}(t_{o})) \, = \, 0$ and $\bold G_{\epsilon}(t_{o}) \; = \; \Theta$, the zero element of $\mathcal{H}$. Thus $supp(\bold G_{\epsilon}) \, \subseteq \, C^{*}$ since $t_{o}$ was any point of $\mathbb{R}^{n} \setminus C^{*} \; = \; C_{*}$.

For any compact subcone $C' \subset C$ and any $r > 0$ (15) yields

\begin{equation}
\mathcal{N}(\bold G_{\epsilon}(t)) \leq M''(C',r,\epsilon)e^{2 \pi <y,t>}, \; t \in \mathbb{R}^{n}, \; y \in (C' \setminus (C' \cap N(\overline{0},r))),
\end{equation}

\noindent as in (17). Also from (15) $\bold g_{\epsilon}(x+iy) \in L^{1}(\mathbb R^n,\mathcal{H}) \cap L^{2}(\mathbb R^n,\mathcal{H})$ as a function of $x \in \mathbb{R}^{n}$ for $y \in C$. Thus from (16) $e^{-2 \pi <y,\cdot>}\bold G_{\epsilon} = \mathcal{F}^{-1}[\bold g_{\epsilon}(x+iy);\cdot], \; y \in C; \; e^{-2 \pi <y,\cdot>}\bold G_{\epsilon} \, \in \, L^{2}(\mathbb R^n,\mathcal{H}), \; y \in C$; and in $L^{2}(\mathbb R^n,\mathcal{H})$

\begin{equation}
\bold g_{\epsilon}(x+iy) = \mathcal{F}[e^{-2 \pi <y,t>}\bold G_{\epsilon}(t);x], \; z = x+iy \; \in \; T^{C}.
\end{equation}

\noindent Now $\bold G_{\epsilon}$ is continuous, $supp(\bold G_{\epsilon}) \subseteq C^{*}$, and (18) holds; thus by Lemma 3.1 $e^{-2 \pi <y,\cdot>}\bold G_{\epsilon} \; \in \; L^{p}(\mathbb R^n,\mathcal{H})$ for all $p, 1 \leq p < \infty$, and for $y \in C$. Thus the Fourier transform in (19) can be interpreted in the $L^{1}(\mathbb R^n,\mathcal{H})$ sense as well as in the $L^{2}(\mathbb R^n,\mathcal{H})$ sense, and (19) becomes

\begin{equation}
\bold g_{\epsilon}(x+iy) = \int_{\mathbb{R}^{n}} \bold G_{\epsilon}(t) e^{2 \pi i<x+iy,t>} dt, \; z = x+iy \in T^{C}.
\end{equation}

Both $\bold G_{\epsilon}$ and $e^{-2 \pi <y,\cdot>}\bold G_{\epsilon}, \; y \in C$, are elements of $\mathcal{S}'(\mathbb R^n,\mathcal{H})$. Also $\bold g_{\epsilon}(\cdot+iy) \in \mathcal{S}'(\mathbb R^n,\mathcal{H}), \; y \; \in \; C$. Let $\phi \in \mathcal{S}$ and $\psi = \mathcal{F}[\phi(t);\cdot] $. We have

\begin{eqnarray}
&~& \langle \bold g_{\epsilon}(x+iy), \psi (x) \rangle = \langle e^{-2 \pi <y,t>}\bold G_{\epsilon}(t), \phi (t) \rangle \\ &~& \rightarrow \langle \bold G_{\epsilon}(t), \phi (t) \rangle = \langle \mathcal{F}[\bold G_{\epsilon}], \psi \rangle \nonumber
\end{eqnarray}

\noindent as $y \; \rightarrow \; \overline{0}, \; y \in C$. As noted previously in this proof, by Theorem 4.1 there is a unique $\bold U_{\epsilon} \in \mathcal{S}'(\mathbb R^n,\mathcal{H})$ such that $\bold g_{\epsilon}(\cdot+iy) \; \rightarrow \; \bold U_{\epsilon}$ in $\mathcal{S}'(\mathbb R^n,\mathcal{H})$ as $y \; \rightarrow \; \overline{0}, \; y \in C$; hence $\mathcal{F}^{-1}[\bold U_{\epsilon}] \; \in \; \mathcal{S}'(\mathbb R^n,\mathcal{H}) $. Thus, $\bold G_{\epsilon} \; = \; \mathcal{F}^{-1}[\bold U_{\epsilon}] \; \in \; \mathcal{S}'(\mathbb R^n,\mathcal{H})$. Moreover, in the sense of the convergence in that space,

\begin{equation}
\lim_{y \rightarrow \overline{0}, y \in C} \bold g_{\epsilon}(x+iy) \; = \; \bold U_{\epsilon} \; = \; \mathcal{F}[\bold G_{\epsilon}] \; \in \mathcal{S}'(\mathbb R^n,\mathcal{H}).
\end{equation}

Recalling the definition of $\bold g_{\epsilon}(z), \; z \in T^{C}$, we have $\bold f(z) \, = \, \bold g_{\epsilon}(z)X_{\epsilon}(z), \, z \in T^{C}, \, \epsilon > 0$, and $\bold f(\cdot+iy)$ has boundary value $\bold h \, \in \, L^{2}(\mathbb R^n,\mathcal{H})$ in $\mathcal{S}'(\mathbb R^n,\mathcal{H})$ as $y \, \rightarrow \, \overline{0}, \; y \, \in \, C$, and $X_{\epsilon}(\cdot+iy)\bold g_{\epsilon}(\cdot+iy) \; \rightarrow \; X_{\epsilon} \mathcal{F}[\bold G_{\epsilon}] \; = \; X_{\epsilon}\bold U_{\epsilon} $ in $\mathcal{S}'(\mathbb R^n,\mathcal{H})$ as $y \, \rightarrow \, \overline{0}, \; y \, \in \, C$. Thus $X_{\epsilon}\bold U_{\epsilon} \; = \; \bold h$ in $\mathcal{S}'(\mathbb R^n,\mathcal{H}), \, \epsilon > 0$. Now for $\phi \, \in \, \mathcal{S}$

\begin{center}
$\langle \frac{\bold h(x)}{X_{\epsilon}(x)}, \phi(x) \rangle = \langle \bold h(x), \frac{\phi(x)}{X_{\epsilon}(x)} \rangle =  \langle \bold U_{\epsilon}, \phi \rangle$
\end{center}
and
$\bold U_{\epsilon} = \frac{\bold h(x)}{X_{\epsilon}(x)} \; \in \; \mathcal{S}'(\mathbb R^n,\mathcal{H}).
$
Clearly,   $\bold H_{\epsilon} = \mathcal{F}^{-1}[\bold h(x)/X_{\epsilon}(x);\cdot] \; \in \; L^{2}(\mathbb R^n,\mathcal{H})$. Since $supp(\bold G_{\epsilon}) \subseteq C^{*}$ then $supp(\bold H_{\epsilon}) \subseteq C^{*}$ almost everywhere. For the function $d_{y}(t)$ defined in (2) we have $d_{y}(t)e^{2 \pi i<z,t>} \; \in \; \mathcal{S}, \; z \in T^{C}$. Thus

\begin{eqnarray}
&~& \int_{C^{*}} \bold G_{\epsilon}(t) e^{2 \pi i<z,t>} dt = \langle \bold G_{\epsilon}(t), d_{y}(t) e^{2 \pi i<z,t>} \rangle \\ &~& = \langle \bold U_{\epsilon}, \mathcal{F}^{-1}[d_{y}(t)e^{2 \pi i<z,t>};\eta] \rangle \nonumber \\ &~& = \langle \frac{\bold h(\eta)}{X_{\epsilon}(\eta)}, \mathcal{F}^{-1}[d_{y}(t)e^{2 \pi i<z,t>};\eta] \rangle \nonumber \\ &~& = \langle \bold H_{\epsilon}(t), d_{y}(t) e^{2 \pi i<z,t>} \rangle = \int_{C^{*}} \bold H_{\epsilon}(t) e^{2 \pi i<z,t>} dt \nonumber
\end{eqnarray}

\noindent with $\bold H_{\epsilon} \in L^{2}(\mathbb R^n,\mathcal{H})$ and $z \in T^{C}$. From [5, Lemma 2.1] $I_{C^{*}}(t)e^{2 \pi i<z,t>} \in  L^{p}$ for all $p, \; 1 \leq p \leq \infty, \; z \in T^{C}$, where $I_{C^{*}}(t)$ is the characteristic function of $C^{*}$. Recalling (20) and (23) we have for $z \in T^{C}$

\begin{eqnarray}
&~& \bold g_{\epsilon}(x+iy) = \int_{C^{*}} \bold G_{\epsilon}(t) e^{2 \pi i<z,t>}dt \\ &~& = \int_{C^{*}} \bold H_{\epsilon}(t) e^{2 \pi i<z,t>}dt = \langle \bold H_{\epsilon}(t), I_{C^{*}}(t) e^{2 \pi i<z,t>} \rangle \nonumber \\ &~& = \langle \mathcal{F}^{-1}[\bold h(x)/X_{\epsilon}(x);t],I_{C^{*}}(t)e^{2 \pi i<z,t>} \rangle \nonumber \\ &~& = \langle \bold h(\eta)/X_{\epsilon}(\eta), \mathcal{F}^{-1}[I_{C^{*}}(t)e^{2 \pi i<z,t>};\eta] \rangle \nonumber \\ &~& = \langle \bold h(\eta)/X_{\epsilon}(\eta), \int_{C^{*}} e^{2 \pi i<z-\eta,t>}dt \rangle= \int_{\mathbb{R}^{n}} \frac{\bold h(t)}{X_{\epsilon}(t)} K(z-t)dt. \nonumber
\end{eqnarray}

Now let $w$ be an arbitrary point of $T^{C}$ and for this arbitrary $w$ consider the function $K(z+w) \bold g_{\epsilon}(z), \; z \in T^{C}$. Using [5, Lemma 3.2] we have $K(z+w)\bold g_{\epsilon}(z)$ is analytic in $T^{C}, \; |K(z+w)| \leq M_{Im(w)} < \infty, \; z \in T^{C}$, where $M_{Im(w)}$ is a constant which depends only on $Im(w)$. Thus $K(z+w)\bold g_{\epsilon}(z) = K(z+w)\bold f(z)/X_{\epsilon}(z)$ satisfies the growth of $\bold f(z), \; z \in T^{C}$, and

\begin{center}
$ \lim_{y \rightarrow \overline{0}, y \in C} K(x+iy+w)\bold g_{\epsilon}(x+iy) = K(x+w) \bold U_{\epsilon} = \frac{K(x+w)\bold h(x)}{X_{\epsilon}(x)}$
\end{center}

\noindent in $\mathcal{S}'(\mathbb R^n,\mathcal{H})$ with $K(x+w)\bold h(x)/X_{\epsilon}(x) \in L^{2}(\mathbb R^n,\mathcal{H})$ since both $K(x+w)$ and $1/X_{\epsilon}(x)$ are bounded for $x \in \mathbb{R}^{n}$. Combining the facts in this paragraph, the same proof leading to (24) applied to $K(z+w)\bold g_{\epsilon}(z),
\; z \in T^{C}$, yields

\begin{equation}
K(z+w)\bold g_{\epsilon}(z) = \int_{\mathbb{R}^{n}} \frac{\bold h(t)}{X_{\epsilon}(t)} K(t+w)K(z-t)dt, \; z \; \in \; T^{C}.
\end{equation}

\noindent Recalling that $w$ above is an arbitrary point of $T^{C}$ we choose $w = -x+iy \in T^{C}$ for $z = x+iy \in T^{C}$. Then (25) combined with (24) becomes

\begin{eqnarray}
&~& \bold g_{\epsilon}(z) = \int_{C^{*}} \bold G_{\epsilon}(t)e^{2 \pi i<z,t>}dt \\ &~& = \int_{\mathbb{R}^{n}} \frac{\bold h(t)}{X_{\epsilon}(t)} K(z-t)dt = \int_{\mathbb{R}^{n}} \frac{\bold h(t)}{X_{\epsilon}(t)} Q(z;t)dt, \; z \in T^{C}. \nonumber
\end{eqnarray}

We now need to construct the function from which the conclusion of this theorem will follow. First we need to show that $\bold h/X_{\epsilon} \; \rightarrow \; \bold h$ in $L^{2}(\mathbb R^n,\mathcal{H})$ as $\epsilon \; \rightarrow \; 0$. Since $|1/X_{\epsilon}(x)| \leq 1, \; x \in \mathbb{R}^{n}, \; \epsilon > 0$, note that

\begin{eqnarray}
&~& (\mathcal{N}(\bold h(x)/X_{\epsilon}(x) \, - \, \bold h(x)))^{2} \; \leq \; (\mathcal{N}(\bold h(x)/X_{\epsilon}(x)) \, + \, \mathcal{N}(\bold h(x)))^{2} \nonumber \\ &~& \leq 2^{2} ((\mathcal{N}(\bold h(x)/X_{\epsilon}(x)))^{2} \, + \, (\mathcal{N}(\bold h(x)))^{2}) \, \leq \, 4((\mathcal{N}(\bold h(x)))^{2} \, + \, (\mathcal{N}(\bold h(x)))^{2}) \nonumber \\ &~& = 8(\mathcal{N}(\bold h(x)))^{2} \nonumber
\end{eqnarray}

\noindent and the right side is independent of $\epsilon > 0$. Also

\begin{center}
$\lim_{\epsilon \rightarrow 0+} \mathcal{N}(\frac{\bold h(x)}{X_{\epsilon}(x)} - \bold h(x)) = \lim_{\epsilon \rightarrow 0+} |(X_{\epsilon}(x))^{-1} - 1| \mathcal{N}(\bold h(x)) \; = \; 0, \; x \in \mathbb{R}^{n}.$
\end{center}

\noindent By the Lebesgue dominated convergence theorem

\begin{equation}
\lim_{\epsilon \rightarrow 0+} |\frac{\bold h(x)}{X_{\epsilon}(x)} - \bold h(x)|_{2} = \lim_{\epsilon \rightarrow 0+} (\int_{\mathbb{R}^{n}}(\mathcal{N}(\frac{\bold h(x)}{X_{\epsilon}(x)} - \bold h(x)))^{2}dx)^{1/2} \; = \; 0.
\end{equation}

Now put

\begin{equation}
\bold G(z) \; = \; \int_{\mathbb{R}^{n}} \bold h(t)Q(z;t)dt, \; z \in T^{C}.
\end{equation}

\noindent Let $z_{o}$ be an arbitrary but fixed point of $T^{C}$. Choose the closed neighborhood $N(z_{o},\rho) \, = \, \{z: |z-z_{o}| \leq \rho, \: \rho > 0 \} \subset T^{C}$ of Lemma 3.3. Using (26), (28), H$\ddot{o}$lder's inequality, and Lemma 3.3

\begin{eqnarray}
&~& \mathcal{N}(\bold g_{\epsilon}(z) - \bold G(z)) \leq (\int_{\mathbb{R}^{n}} (\mathcal{N}(\frac{\bold h(t)}{X_{\epsilon}(t)} - \bold h(t)))^{2}dt)^{1/2} ||Q(z;t)||_{L^{2}} \\ &~& \leq B(z_{o}) (\int_{\mathbb{R}^{n}} (\mathcal{N}(\frac{\bold h(t)}{X_{\epsilon}(t)} - \bold h(t)))^{2}dt)^{1/2} \nonumber
\end{eqnarray}

\noindent for $z \in N(z_{o},\rho) \subset T^{C}$. (27) and (29) now yield

\begin{center}
$\lim_{\epsilon \rightarrow 0+} \; \bold g_{\epsilon}(z) \; = \; \bold G(z)$
\end{center}

\noindent uniformly in $z \in N(z_{o},\rho)$. Since $\bold g_{\epsilon}(z)$ is analytic in $T^{C}$ for each $\epsilon > 0$, we have that $\bold G(z)$ is analytic at $z_{o} \in T^{C}$ and hence in $T^{C}$ since $z_{o}$ is an arbitrary point of $T^{C}$. Applying Lemma 3.4 we have $\bold G(z)$ of (28) is an element of $H^{2}(T^{C},\mathcal{H})$.

For $\phi \in \mathcal{S}$ we use H$\ddot{o}$lder's inequality to obtain

\begin{eqnarray}
&~& \mathcal{N}( \langle \bold G(x+iy),\phi(x) \rangle - \langle \bold h(x),\phi(x) \rangle ) = \mathcal{N}(\int_{\mathbb{R}^{n}}(\bold G(x+iy) - \bold h(x)) \phi(x) dx) \nonumber \\ &~& \leq |\bold G(x+iy) - \bold h(x)|_{2} ||\phi||_{L^{2}}.
\end{eqnarray}

\noindent By Lemma 3.4 $\bold G(x+iy) \rightarrow \bold h(x)$ in $L^{2}(\mathbb R^n,\mathcal{H})$ as $y \rightarrow \overline{0}, \; y \in C$; hence  $\bold G(x+iy) \rightarrow \bold h(x)$ in $\mathcal{S}'(\mathbb R^n,\mathcal{H})$ as $y \rightarrow \overline{0}, \; y \in C$.

We now consider $\bold f(z) - \bold G(z), \; z \in T^{C}$, which is analytic in $T^{C}$; $f$ satisfies the growth (12) and $\bold G$ satisfies the growth (6). Thus

\begin{eqnarray}
&~& \mathcal{N}(\bold f(z) - \bold G(z)) \leq P(C',r)(1+|z|)^{R}, \\ &~& z = x+iy \in T(C',r) = \mathbb{R}^{n}+i(C' \setminus (C' \cap N(\overline{0},r))), \nonumber
\end{eqnarray}

\noindent for any compact subcone $C' \subset C$ and any $r > 0$, where $P(C',r)$ is a constant depending on $C' \subset C$ and on $r > 0$, and

\begin{equation}
\lim_{y \rightarrow \overline{0}} \; (\bold f(x+iy) - \bold G(x+iy)) \; = \; \bold h(x) - \bold h(x) \; = \; \Theta
\end{equation}

\noindent in $\mathcal{S}'(\mathbb R^n,\mathcal{H}))$. Now put $\bold F(z) = \bold f(z)-\bold G(z), \; z \in T^{C}$, and $\bold F(z)$ satisfies (31) and (32). Letting $\epsilon = 1$ in the function $X_{\epsilon}(z)$ at the beginning of this proof, consider $\bold g(z) = \bold F(z)/X_{1}(z), z \in T^{C}.$ As in (15), for any compact subcone $C' \subset C$ and any $r > 0$

\begin{eqnarray}
&~& \mathcal{N}(\bold g(z) \leq P'(C',r)(1+|z|)^{-n-2}, \nonumber \\ &~& z = x+iy \in T(C',r) = \mathbb{R}^{n}+i(C' \setminus (C' \cap N(\overline{0},r))). \nonumber
\end{eqnarray}

\noindent Now putting as in (16)

\begin{center}
$\bold  A(t) \; = \; \int_{\mathbb{R}^{n}} \bold g(x+iy) e^{-2 \pi i<x+iy,t>}dx, \; y \in C, \; t \in \mathbb{R}^{n},$
\end{center}

\noindent and proceeding with the form of the proof from (16) to (20) we have that $\bold A$ is continuous, is independent of $y \in C$, has support in $C^{*}$, satisfies a growth as in (18), $e^{-2 \pi <y,t>}\bold A(t) \, = \, \mathcal{F}^{-1}[\bold g(x+iy);t], \; t\in\mathbb R^n, \; y \in C$, with $e^{-2 \pi <y,\cdot>}\bold A \in L^{2}(\mathbb R^n,\mathcal{H}), \; \bold g(x+iy) = \mathcal{F}[e^{-2 \pi <y,t>}\bold A(t);x], \; x\in\mathbb R^n, \; y \in C$, and

\begin{equation}
\bold g(x+iy) = \int_{\mathbb{R}^{n}} \bold A(t) e^{2 \pi i<x+iy,t>} dt, \; z = x+iy \in T^{C}.
\end{equation}

\noindent For $\phi \in \mathcal{S}$ and $y \in C$

\begin{eqnarray}
&~& \langle \bold F(z)/X_{1}(z),\phi(x) \rangle = \langle \bold g(z), \phi(x) \rangle \nonumber \\ &~& = \langle \mathcal{F}[e^{-2 \pi <y,t>}\bold A(t);x],\phi(x) \rangle = \langle e^{-2 \pi <y,t>}\bold A(t),\hat{\phi}(t) \rangle \nonumber \\ &~& \rightarrow \langle \bold A(t),\hat{\phi}(t) \rangle = \langle \mathcal{F}[\bold A],\phi \rangle \nonumber
\end{eqnarray}

\noindent as $y \rightarrow \overline{0}, \; y \in C$. Thus

\begin{eqnarray}
&~& \langle \bold F(x+iy),\phi(x) \rangle = \langle \bold g(x+iy),X_{1}(x+iy) \phi(x) \rangle \nonumber \\ &~& \rightarrow \langle \mathcal{F}[\bold A],X_{1}(x) \phi(x) \rangle = \langle X_{1}(x) \mathcal{F}[\bold A], \phi (x) \rangle \nonumber
\end{eqnarray}

\noindent as $y \rightarrow \overline{0}, \; y \in C$. Combining this fact with (32) we have $X_{1}(x) \mathcal{F}[\bold A] = \Theta$  which yields $\bold A = \Theta$ in $\mathcal{S}'(\mathbb R^n,\mathcal{H})$.

Now put

\begin{center}
$\Delta \; = \; \prod_{j=1}^{n} \left( 1-i(-1)^{v_{j}} \left( \frac{-1}{2 \pi i} \frac{\partial}{\partial t_{j}} \right) \right)^{R+n+2}.$
\end{center}

\noindent From (33) and for $z \in T^{C}$

\begin{eqnarray}
&~& \bold F(z) = X_{1}(z) \int_{\mathbb{R}^{n}} \bold A(t)e^{2 \pi i<z,t>}dt \nonumber \\ &~& = \langle \Delta \bold A(t), d_{y}(t)e^{2 \pi i<z,t>} \rangle = \bold \Theta \nonumber
\end {eqnarray}

\noindent which yields

\begin{center}
$\bold f(z) = \bold G(z) = \int_{\mathbb{R}^{n}} \bold h(t)Q(z;t)dt, \; z \in T^{C},$
\end{center}

\noindent and $\bold f(z) \in H^{2}(T^{C},\mathcal{H})$. The proof of Theorem 4.2 is complete. \\

Extending Theorem 4.2 to the cases $2 \, < \, p \, \leq \, \infty$ we have the following result. \\

THEOREM 4.3. {\textit{Let $C$ be an open convex cone which is contained in or is any of the $2^{n}$ n-rants $C_{v}$ in $\mathbb{R}^{n}$. Let $\bold f(z)\in
\mathcal A(T^{C},\mathcal{H})$  which satisfies (12). Let the unique boundary value $\bold U$ of Theorem 4.1 be $\bold h \, \in \, L^{p}(\mathbb R^n,\mathcal{H}), \; 2 \, < \, p \, \leq \, \infty$. We have $\bold f \; \in \; H^{p}(T^{C},\mathcal{H}), \; 2 < p \leq \infty$, and}}

\begin{center}
$ \bold f(z) = \int_{\mathbb{R}^{n}} \bold h(t)Q(z;t)dt, \; z \in T^{C}.$ \\
\end{center}

{\textit{Proof.}} Note that the analysis from the beginning of the proof of Theorem 4.2 through the fact that $\bold U_{\epsilon} = \frac{\bold h(x)}{X_{\epsilon}(x)} \in \mathcal{S}'(\mathbb{R}^{n},\mathcal{H})$ above (23) is independent of the value of $p, \, 2 \, \leq \, p \, \leq \, \infty$. From the definition of $X_{\epsilon}(z)$ at the beginning of the proof of Theorem 4.2

\begin{center}
$|1/X_{\epsilon}(x)| = \epsilon^{-n(R+n+2)} \prod_{j=1}^{n} (\epsilon^{-2}+x_{j}^{2})^{-1-R/2-n/2}, \; x \in \mathbb{R}^{n},$
\end{center}

\noindent and $1/X_{\epsilon} \in L^{q}, \; 1 \leq q \leq \infty$. Thus $\bold h/X_{\epsilon} \in L^{1}(\mathbb R^n,\mathcal{H}) \cap L^{p}(\mathbb R^n,\mathcal{H}), \; 2 < p \leq \infty$. If $p = \infty, \bold h(x)/X_{\epsilon}(x) \in L^{1}(\mathbb{R}^{n}, \mathcal{H}) \cap L^{2}(\mathbb{R}^{n},\mathcal{H}) \cap L^{\infty}(\mathbb{R}^{n},\mathcal{H})$. Further, if $2 < p < \infty$

\begin{eqnarray}
&~& \int_{\mathbb{R}^{n}} (\mathcal{N}(\bold h(x)/X_{\epsilon}(x)))^{2}dx = \int_{\mathbb{R}^{n}} (\mathcal{N}(\bold h(x))^{2}|1/X_{\epsilon}(x)|^{2}dx \nonumber \\ &~& \leq ||(\mathcal{N}(\bold h(x)))^{2}||_{L^{p/2}} || \; |1/X_{\epsilon}(x)|^{2} \; ||_{L^{p/(p-2)}} \; < \; \infty. \nonumber
\end{eqnarray}

\noindent Thus $\bold h/X_{\epsilon} \in L^{1}(\mathbb{R}^{n},\mathcal{H}) \cap L^{2}(\mathbb R^n,\mathcal{H}) \cap L^{p}(\mathbb R^n,\mathcal{H})$ for the value of $p, 2 < p \leq \infty$. As noted before in this proof the analysis from the beginning of the proof of Theorem 4.2 through the fact that $\bold U_{\epsilon} = \frac{\bold h(x)}{X_{\epsilon}(x)} \in \mathcal{S}'(\mathbb{R}^{n},\mathcal{H})$ above (23) is independent of the value of $p$ in the hypotheses of Theorem 4.2 and also now of the hypotheses of Theorem 4.3. Since $\bold h/X_{\epsilon} \in L^{2}(\mathbb R^n,\mathcal{H})$ for $\bold h \in L^{p}(\mathbb R^n,\mathcal{H}), \; 2  <  p  \leq  \infty$, here, the analysis from this fact that $\bold U_{\epsilon} = \frac{\bold h(x)}{X_{\epsilon}(x)} \in \mathcal{S}'(\mathbb{R}^{n},\mathcal{H})$ to equation (23) to equation (26) holds exactly as in the proof of Theorem 4.2 for the present case. For $2  <  p  <  \infty$

\begin{eqnarray}
&~& (\mathcal{N}(\bold h(x)/X_{\epsilon}(x) - \bold h(x)))^{p} \leq (\mathcal{N}(\bold h(x)/X_{\epsilon}(x)) + \mathcal{N}(\bold h(x)))^{p} \nonumber \\ &~& \leq 2^{p}((\mathcal{N}(\bold h(x)/X_{\epsilon}(x)))^{p} + (\mathcal{N}(\bold h(x)))^{p}) \leq 2^{p+1}(\mathcal{N}(\bold h(x)))^{p}. \nonumber
\end{eqnarray}

\noindent The use of the Lebesgue dominated convergence theorem as in (27) now shows

\begin{equation}
\lim_{\epsilon \rightarrow 0+} |\bold h(x)/X_{\epsilon}(x) - \bold h(x)|_{p} \; = \; 0.
\end{equation}

\noindent Define $\bold G$ as in (28) and recall that Lemma 3.3 holds for all $p, \; 1 \leq p \leq \infty$. For the present case of $2 < p < \infty$ we use an estimate as in (29) together with (34) to obtain that $\bold G$ is analytic in $T^{C}$. For the case $p = \infty$ choose a closed neighborhood contained in $T^{C}$ about each fixed point $z_{o} \in T^{C}$, which can be done since $C$ is open; and $\bold h/X_{\epsilon} \in L^{2}(\mathbb{R}^{n},\mathcal{H})$ for $\bold h \in L^{\infty}(\mathbb{R}^{n},\mathcal{H})$. From (26)

\begin{center}
$\bold g_{\epsilon}(z) = \int_{\mathbb{R}^{n}} \frac{\bold h(t)}{X_{\epsilon}(t)}Q(z;t)dt, \; z \in T^{C},$
\end{center}

\noindent is analytic in $T^{C}$, and from Lemma 3.5 $\bold g_{\epsilon}(z) \rightarrow \bold G(z)$ uniformly on the closed neighborhood about $z_{o} \in T^{C}$ and hence on the corresponding open neighborhood about $z_{o} \in T^{C}$ as $\epsilon \rightarrow 0+$. Thus for the case $\bold h \in L^{\infty}(\mathbb{R}^{n},\mathcal{H})$, $\bold G$ is analytic on $T^{C}$ since $z_{o}$ was an arbitrary point of $T^{C}$. Since $\bold h \in L^{p}(\mathbb{R}^{n},\mathcal{H}), 2 < p \leq \infty$, here, Lemma 3.4 yields $\bold G \in H^{p}(T^{C},\mathcal{H}), 2 < p \leq \infty$. Arguing as in (30) we have for $2 < p < \infty, 1/p \; + \; 1/q \; = \; 1,$

\begin{center}
$\mathcal{N}( \langle \bold G(x+iy), \phi(x) \rangle - \langle \bold h(x), \phi(x) \rangle ) \leq |\bold G(x+iy)-\bold h(x)|_{p} ||\phi||_{L^{q}}, \; \phi \in \mathcal{S},$
\end{center}

\noindent and using Lemma 3.4 we get $\bold G(x+iy) \; \rightarrow \; \bold h(x)$ in $\mathcal{S}'(\mathbb R^n,\mathcal{H})$ as $y \; \rightarrow \; \overline{0},\; y \in C$; and this convergence holds also in $\mathcal{S}'(\mathbb{R}^{n},\mathcal{H})$ for $p = \infty$.

We now have both $\bold f$ and $\bold G$ are analytic in $T^{C}$, and $\bold G(z)$ satisfies the stated growth in Lemma 3.4 for $z \in T(C',r)$ in the cases $2 < p < \infty$ and is bounded for all $z \in T^{C}$ independent of $z$ if $p = \infty$. Further both $\bold f$ and $\bold G$ have $\bold h \in L^{p}(\mathbb R^n,\mathcal{H})$ as boundary value in $\mathcal{S}'(\mathbb R^n,\mathcal{H})$. Thus we may consider $(\bold f(z) - \bold G(z)), \; z \in T^{C}$, exactly as in the proof of Theorem 4.2 starting at (31) and continuing through the end of the proof for the case $p = 2$ to conclude that

\begin{center}
$\bold f(z) = \bold G(z) = \int_{\mathbb{R}^{n}} \bold h(t)Q(z;t)dt, \; z \in T^{C},$
\end{center}

\noindent and $\bold f \in H^{p}(T^{C}, \mathcal{H}), \; 2 < p \leq \infty$, from Lemma 3.4. The proof of Theorem 4.3 is complete. \\

We now extend Theorems 4.2 and 4.3 to a tube $T^{C}$ where $C$ is an arbitrary regular cone. \\

THEOREM 4.4. {\textit{Let $C$ be a regular cone. Let $\bold f\in\mathcal A(T^{C},\mathcal{H})$ and satisfy (12). Let the unique $\mathcal{S}^{(m)'}(\mathbb R^n,\mathcal{H}) \subset \mathcal{S}'(\mathbb R^n,\mathcal{H})$ boundary value of $\bold f$ from Theorem 4.1 be $\bold h \in L^{p}(\mathbb R^n,\mathcal{H}), \; 2 \leq p \leq \infty$. We have $\bold f \in H^{p}(T^{C}, \mathcal{H}), \; 2 \leq p \leq \infty$.}} \\

{\textit{Proof.}} For each of the $2^{n}$ n-rants $C_{v}$ consider $C \cap C_{v}$. Let $S_{j}, \; j = 1,...,k$, be an enumeration of the  intersections $C \cap C_{v}$ which are non-empty; each $S_{j}$ is an open regular cone which is contained in or is a n-rant $C_{v}$ in $\mathbb{R}^{n}$. Put

\begin{equation}
\bold f_{j}(z) = \bold f(z), \; z \in T^{S_{j}} = \mathbb{R}^{n}+iS_{j}, \; j = 1,...,k.
\end{equation}

\noindent Each $\bold f_{j}$ satisfies the analyticity and growth hypotheses of Theorems 4.2 and 4.3, and each $\bold f_{j}$ obtains the unique $\bold h \in L^{p}(\mathbb R^n,\mathcal{H})$ as $\mathcal{S}'(\mathbb R^n,\mathcal{H})$ boundary value. By Theorems 4.2 and 4.3 each $\bold f_{j} \in H^{p}(T^{S_{j}}, \mathcal{H}), \; 2 \leq p \leq \infty, \; j = 1,...,k$; and

\begin{equation}
f\bold (z) = \bold f_{j}(z) = \int_{\mathbb{R}^{n}} \bold h(t)Q(z;t)dt, \; z \in T^{S_{j}}, \; j = 1,...,k.
\end{equation}

For $2 \leq p < \infty$ there are constants $A_{j}, \; j = 1,...,k$, independent of $y = Im(z)$ such that

\begin{equation}
\int_{\mathbb{R}^{n}} (\mathcal{N}(\bold f(x+iy)))^{p}dx = \int_{\mathbb{R}^{n}} (\mathcal{N}(\bold f_{j}(x+iy)))^{p}dx \leq A^{p}_{j} \; , \; y \in S_{j}, \; j = 1,...,k.
\end{equation}

\noindent Put

\begin{equation}
A \; = \; max \{ A_{1},A_{2},...,A_{k} \}.
\end{equation}

\noindent Now let $y \in C$ such that $y \notin S_{j}, j = 1,...,k$. Then $y$ is on the topological boundary of some $S_{j}$. For this $S_{j}$ choose a sequence of points $\{ y_{j,l} \} \subset S_{j}$ such that $y_{j,l} \rightarrow y$ as $l \rightarrow \infty$. Since $\bold f(z)$ is analytic in $T^{C}$, by Fatou's lemma we have for $y \in C$ such that $y \notin S_{j}, \; j = 1,...,k$,

\begin{eqnarray}
&~& \int_{\mathbb{R}^{n}} (\mathcal{N}(\bold f(x+iy)))^{p}dx \leq \liminf_{l \rightarrow \infty} \int_{\mathbb{R}^{n}} (\mathcal{N}(\bold f(x+iy_{j,l})))^{p}dx \nonumber \\ &~& \leq A^{p}_{j} \; \leq \; A^{p}.
\end{eqnarray}

\noindent Combining (35), (37), (38), and (39) we have for $2 \leq p < \infty$

\begin{center}
$ \int_{\mathbb{R}^{n}} (\mathcal{N}(\bold f(x+iy)))^{p}dx \; \leq \; A^{p}, \; y \in C$,
\end{center}

\noindent where $A$ is independent of $y \in C$. Thus $\bold f \in H^{p}(T^{C},\mathcal{H}), 2 \leq p < \infty.$

For $p = \infty$ each $\bold f_{j}(z) \in H^{\infty}(T^{S_{j}},\mathcal{H}), j = 1,...,k$; and $\mathcal{N}(\bold f_{j}(z)) \leq B_{j}, z \in T^{S_{j}}, j = 1,...,k,$ for positive constants $B_{j}, j = 1,...,k$, which are independent of $z \in T^{S_{j}}, j = 1,...,k$. Put $B = max \{B_{1},B_{2},...,B_{k} \}$. If $y \in C$ such that $y \notin S_{j}, j = 1,...,k,$ again choose a sequence $\{ y_{j,l} \} \subset S_{j}$ for some $j$ such that $y_{j,l} \rightarrow y$ as $l \rightarrow \infty$; since $\bold f(z)$ is analytic and hence continuous, a simple continuity argument yields $\mathcal{N}(\bold f(x+iy)) \; \leq \; 1+B$ for any $y \in C$ such that $y \notin S_{j}, j = 1,...,k$. Combining these facts we thus have $\mathcal{N}(\bold f(x+iy)) \; \leq \; 1+B$ for all $z \in T^{C}$, and $\bold f \in H^{\infty}(T^{C},\mathcal{H})$.

We conclude $\bold f \in H^{p}(T^{C},\mathcal{H})$ for each choice of $p, \; 2 \leq p \leq \infty$, and the proof of Theorem 4.4 is complete. \\

\noindent {\textbf{5 Results for $1 \leq p < 2$}} \\

For the cases $1 \leq p < 2$ Theorem 4.3, and hence Theorem 4.4, will not follow from the $p = 2$ case of Theorem 4.2 by a proof like the one used to obtain the $2 < p \leq \infty$ cases from the $p = 2$ case in this paper. To obtain results like Theorems 4.2, 4.3, and 4.4 for the cases $1 \leq p < 2$ we must use separate proof techniques. We propose to do so by using the concept of Banach space $\mathcal{B}$ with Fourier type $p, 1 \leq p < 2$, as discussed in [14, section 6].

From [14, section 6] the Banach space $\mathcal{B}$ has Fourier type $p$ with respect to $\mathbb{R}^{n}$ if there is a constant $K > 0$ such that for every compactly supported $\bold f \in L^{p}(\mathbb{R}^{n},\mathcal{B}), 1 \leq p \leq 2$,

\begin{equation}
(\int_{\mathbb{R}^{n}} (\mathcal{N}(\widehat{\bold f}(t)))^{q}dt)^{1/q} \; \leq \; K(\int_{\mathbb{R}^{n}} (\mathcal{N}(\bold f(t)))^{p}dt)^{1/p}, \; q \; = \; p/(p-1).
\end{equation}

\noindent By completeness, for a Banach space with Fourier type $p, 1 \leq p \leq 2$, we have for every $\bold f \in L^{p}(\mathbb{R}^{n},\mathcal{B})$ that $\widehat{\bold f} \in L^{q}(\mathbb{R}^{n},\mathcal{B})$ and (40) holds. Again we refer to [14, section 6] for details.

If $\mathcal{B}$ is of Fourier type $p, 1 \leq p \leq 2$, the Fourier transform of $\bold f$ defined in the sense of the Bochner integral imbedded into $\mathcal{S}'(\mathbb{R}^{n},\mathcal{B})$ coincides with $\bold g = \mathcal{F}[\bold f]$ defined as the distributional Fourier transform, and $\bold f = \mathcal{F}^{-1}[\bold g]$ in $\mathcal{S}'(\mathbb{R}^{n},\mathcal{B})$.

In future research we use the concept of a Banach space $\mathcal{B}$ having Fourier type $p, 1 \leq p \leq 2$, to obtain the results Theorems 4.2, 4.3, and 4.4 for the cases $1 \leq p < 2$. \\

\noindent {\textbf{References}} \\

\noindent [1]  W. Arendt, C.J.K. Batty M. Hieber, F. Neubrander, {\textit{Vector-valued Laplace Transforms and Cauchy Problems, Second edition}},
 Birkhauser/Springer Basel AG, Basel, 2011. \\

\noindent [2]  R.D. Carmichael, W.W. Walker, {\textit{Representation of distributions with compact support}}, Manuscripta Math. 11 (1974), 305 -338. \\

\noindent [3]  R.D. Carmichael, E.K. Hayashi, {\textit{A pointwise growth estimate for analytic functions in tubes}}, Internat. J. Math. Math. Sci. 3 (1980), 575 - 581. \\

\noindent [4]  R.D. Carmichael, S.P. Richters, {\textit{Growth of $H^{p}$ functions in tubes}}, Internat. J. Math. Math. Sci. 4 (1981), 435 - 443. \\

\noindent [5]  R.D. Carmichael, S.P. Richters, {\textit{Holomorphic functions in tubes which have distributional boundary values and which are $H^{p}$ functions}}, SIAM J. Math. Anal. 14 (1983), 596 - 621. \\

\noindent [6]  R.D. Carmichael, D. Mitrovi$\acute{c}$, {\textit{Distributions and Analytic Functions}}, Longman Scientific and Technical, Harlow, UK, 1989. \\

\noindent [7]  R.D. Carmichael, A. Kami$\acute{n}$ski, S. Pilipovi$\acute{c}$, {\textit{Boundary Values and Convolution In Ultradistribution Spaces}}, World Scientific Publishing Co., Singapore, 2007. \\

\noindent [8]  A. Debrouwere, J. Vindas, {\textit{On the non-triviality of certain spaces of analytic functions. Hyperfunctions and ultrahyperfunctions of fast growth}}, Rev. R. Acad. Cienc. Exactas Fis. Nat. Ser. A. Math. RACSAM, in press, doi:10.1007/s13398-017-0392-9. \\

\noindent [9]  G. Debruyne, J. Vindas, {\textit{Complex Tauberian theorems for Laplace transforms with local pseudofunction boundary behavior}}, J. Anal. Math., to appear. \\

\noindent [10] P. Dimovski, S. Pilipovi$\acute{c}$, J. Vindas, {\textit{New distribution spaces associated to translation-invariant Banach spaces}}, Monatsh. Math. 177 (2015), 495 - 515. \\

\noindent [11] P. Dimovski, S. Pilipovi$\acute{c}$, J. Vindas, {\textit{Boundary values of holomorphic functions in translation-invariant distribution spaces}}, Complex Var. Elliptic Equ. 60 (2015), 1169 - 1189. \\

\noindent [12] N. Dunford, J.T. Schwartz, {\textit{Linear Operators Part I: General Theory}}, Interscience Publishers Inc., New York, 1966. \\

\noindent [13] H. Federer, {\textit{Geometric Measure Theory}}, Springer-Verlag, New York, 1969. \\

\noindent [14] J. Garcia-Cuerva, K. S. Kazarian, V. I. Kolyada, J. L. Torrea,
{\textit{Vector-valued Hausdorff-Young inequality and applications}}, Russian Math. Surveys 53 (1998), 435 - 513 (Uspekhi Mat. Nauk 53 (1998), 3 - 84).\\

\noindent [15] R. Meise, {\textit{Darstellung temperierter vektorwertiger Distributionen durch holomorphe Funktionen I}}, Math. Ann. 198 (1972), 147 -159. \\

\noindent [16] R. Meise, {\textit{Darstellung temperierter vektorwertiger Distributionen durch holomorphe Funktionen II}}, Math. Ann. 198 (1972), 161 - 178. \\

\noindent [17] S Pilipovi$\acute{c}$, {\textit{Tempered ultradistributions}}, Boll. Un. Mat. Ital. 7 2-B (1988), 235 - 251. \\

\noindent [18] A.K. Raina, {\textit{On the role of Hardy spaces in form factor bounds}}, Lett. Math. Phys. 2 (1978), 513 - 519. \\

\noindent [19] L. Schwartz, {\textit{Th$\acute{e}$orie des distributions a valeurs vectorielles. I}}, Ann. Inst. Fourier 7 (1957), 1 - 149. \\

\noindent [20] L. Schwartz, {\textit{Th$\acute{e}$orie des distributions a valeurs vectorielles. II}}, Ann. Inst. Fourier 8 (1958), 1 - 209. \\

\noindent [21] E.M. Stein, G. Weiss, {\textit{Introduction to Fourier Analysis on Euclidean Spaces}}, Princeton University Press, Princeton, NJ, 1971. \\

\noindent [22] H.G. Tillmann, {\textit{Darstellung der Schwartzschen Distributionen durch analytische Funktionen}}, Math. Z. 77 (1961), 106 - 124. \\

\noindent [23] V.S. Vladimirov, {\textit{Methods of the Theory of Funtions of Many Complex Variables}}, M.I.T. Press, Cambridge, MA, 1966. \\

\end{document}